\documentclass[a4paper,11pt]{article}
\usepackage[margin=2cm]{geometry}
\usepackage{amsmath}
\usepackage{amsfonts}
\usepackage{amssymb}
\usepackage{graphicx}
\usepackage{dsfont}
\newtheorem{theorem}{Theorem}

\newtheorem{lemma}{Lemma}

\newcommand{\X}{{\cal X}}
\newcommand{\Pb}{{\mathbb{P}}}
\newcommand{\Eb}{{\mathbb{E}}}
\newcommand{\cov}{\text{Cov}}
\title{Consistency of the recursive nonparametric regression estimation for dependent functional data}
 \author{ Aboubacar Amiri\thanks{EQUIPPE} and Baba Thiam
\thanks{ EQUIPPE, Universit\'e Charles De Gaulle, Lille
3, Maison de la Recherche, Domaine universitaire du Pont de Bois, BP 60149, 59653
Villeneuve d'ascq cedex, France.
baba.thiam@univ-lille3.fr}}

\begin{document}

\maketitle

\begin{abstract}
We consider the recursive estimation of a regression functional where the explanatory variables take values in some functional space. We prove the almost sure convergence of such estimates for dependent functional data. Also we derive the mean quadratic error of the  considered class of estimators.  Our results are established with rates and asymptotic appear bounds, under strong mixing condition.\bigskip


\noindent{\bf Classcode:}
62G05, 62G07,  62G08, 62G20, 62L12.\bigskip

\end{abstract}

\section{Introduction}\label{intro}
In this paper we  study  the regression model of a scalar  response variable  given a  functional covariate.
%A functional   random process is a collection  of random variables valued on some functional space.
Functional data analysis  is a problem of considerable interest in statistics and  has been found to be useful in many practical fields, including climatology, economics,
linguistics, medicine,... The statistical study of this kind of data  is the subject of many papers in parametric and nonparametric statistics.  For background material on this subject we highlight  the  works of Ramsay and   Dalzell \cite{RamsayDalzell}, Ramsay and Silverman \cite{RamsaySilverman02,RamsaySilverman05}.  Since these  pioneer   contributions,  the literature on this topic is still growing. A survey of the nonparametric functional regression appears in Ferraty et al. \cite{Ferraty2007}, while more recent results are collected in the book  by Ferraty and Vieu \cite{Ferraty2006}. There are several ways to study  the link between a  response variable  given an explanatory variable.  For example, one of the most studied models is the regression model when the response variable  $Y$ is real and the explanatory variable  ${\cal X}$  belongs to
some functional space ${\cal E}$.  Then, the regression
model writes $Y=r({\cal X})+\varepsilon,$ where $r: {\cal E}\rightarrow \mathbb{R}$ is an operator and $\varepsilon$ is an error random variable. Many
works have been done around this model when the operator $r$  is supposed to be linear, contributing to the popularity of the so-called functional linear
model.
%In this linear context, the operator $r$  writes $\langle \alpha, \cdot\rangle$ where $\langle \cdot, \cdot\rangle$ stands
%for an inner product of the space ${\cal E}$ and $\alpha$
%belongs to ${\cal E}$. The goal is then to estimate the unknown function $\alpha$.
We refer the reader for instance to
the works of Cardot et al. \cite{Cardot03} or Crambes et al. \cite{Crambes09} for different methods to estimate  $r$ in this linear context.
%$\alpha$.
 Another way is to estimate $r$ by a nonparametric approach. The first results on this  context  were obtained by Ferraty and Vieu \cite{Ferraty2000}. They established the almost complete convergence of a  kernel estimator of the regression function in the i.i.d  case. The  study  on their   Nadaraya-Watson type estimator  is extended to several directions.  Dabo-Niang and Rhomari \cite{Dabo03} stated the $L^p$-convergence  of the kernel estimator, while  Delsol \cite{Delsol09} gave the $L^p$-convergence  with asymptotic appear bound. The asymptotic normality of the same estimator has been obtained by Masry \cite{Masry2005} under strong mixing conditions  and extended by Delsol \cite{Delsol09}.  Ling and Wu \cite{Ling2012} stated the almost sure convergence of the kernel estimator under strong mixing conditions.  Functional data appear in many practical situations, as soon as one is interested on a continuous  phenomenon. To consider such  data as objects belonging to some functional space brings  more precisions on  the studied phenomenon. However, the computation of the estimators  can be time consuming in this context, the  use  of recursive methods remains a good alternative to the classical ones. By `recursive', we mean that
the estimator calculated from the first $n$ observations, say $f_n$, is only a function of $f_{n-1}$ and the
$n^{th}$ observation. In this way, the estimator can be updated with each new observation added to
the database. \\  The purpose of this paper is  to  apply recursive methods
to functional data.  Recursive estimation is achieved with the use of recursive estimators, typically kernel ones.  For informations on  nonparametric recursive methods, the reader is referred to the books by Gyorfi et al. \cite{Gyorfi2002}, or the recent works of Vilar and Vilar \cite{Vilar2000},  Wang and Liang \cite{Wang2004},  Quintela-Del-Rio \cite{Quintela2010}, Amiri \cite{Amiri1} and the references there in. The first results concerning  the recursive kernel  estimator of the regression function with  functional explanatory variable were obtained by  Amiri et al. \cite{Amiri2}. They established the mean square error, the almost sure convergence with
rates  and a central limit theorem  for a class of recursive kernel estimates of the regression function when the explanatory variable is functional and the observations are i.i.d.  The main goal of this paper is the extension of a few of the results obtained by  Amiri et al. \cite{Amiri2} to dependent data.  The rest of the paper proceeds as follows.  We will present the regression model on section \ref{model}. On section \ref{resultat},  we give  assumptions and results   on the strong consistency  and mean quadratic error for the recursive regression estimate. Section \ref{preuve} is devoted to the proofs of our results.
\section{Recursive regression  estimate for curves}\label{model}

Let us consider a random process $Z_t=(\X_t,Y_t), t\in \mathbb{N},$ where $Y_t$ is a scalar random variable and $\X_t$ takes values in some functional space $\mathcal{E}$ endowed with a semi-norm $\|\cdot\|$. Assume the existence of an operator $r$ satisfying  $
r(\chi):=\Eb\left(Y_t|\X_t=\chi\right),~~\chi\in \mathcal{E},\text{ for all } t \in \mathbb{N} .
$
%where $$.
%To estimate $r$, one can consider the following kernel recursive estimation
%\begin{equation}
%r_{n} (\chi) := \frac{{\displaystyle \sum_{i=1}^{n} Y_{i} K \left( \frac{\left\| \chi - \mathcal{X}_{i} \right\|}{h_n} \right)}}{{\displaystyle \sum_{i=1}^{n} K \left( \frac{\left\| \chi - \mathcal{X}_{i} \right\|}{h_n} \right)}}, \label{estimFV}
%\end{equation}
%where $K$ is a kernel, $(h_n)$ a sequence of bandwidths. The asymptotic behavior of this estimator proposed by \cite{Ferraty2007} is widely investigated in the literature........\\
To estimate  $r$, one can consider the family of recursive estimators indexed by a parameter $\ell\in[0,1]$ introduced in Amiri et al.  \cite{Amiri2} and defined by
% Nonparametric regression aims to estimate the functional  $r(\chi):=\text{E}\left( Y |{\cal X}=\chi \right),$ for $\chi \in {\cal E}$. To this end, let us consider the family of recursive estimators indexed by a parameter $\ell\in[0,1],$ and defined by
 $$r_n^{[\ell]}({\chi}) :=\frac{\sum\limits_{i=1}^n\frac{Y_i}{F(h_i)^{\ell}}K\left( \frac{\|{\chi}-{\cal X}_i\|}{h_i}\right)}
{\sum\limits_{i=1}^n\frac{1}{F(h_i)^{\ell}}K\left( \frac{\|{\chi}-{\cal X}_i\|}{h_i}\right)},$$
where $K$ is a kernel, $(h_n)$ a sequence of bandwidths and $F$ is the cumulative distribution function of the random variable $\|\chi-{\cal X}\|$.   This family of estimators  is a recursive modification of the Nadaraya-Watson type estimator of Ferraty and Vieu \cite{Ferraty2006} and can be computed recursively by
$$r_{n+1}^{[\ell]}(\chi)=\frac{\left[\sum\limits_{i=1}^nF(h_i)^{1-\ell}\right]\varphi_n^{[\ell]}(\chi)+\left[\sum\limits_{i=1}^{n+1}F(h_i)^{1-\ell}\right]Y_{n+1}K_{n+1}^{[\ell]}\left(\|\chi-{\cal X}_{n+1} \|\right)}{\left[\sum\limits_{i=1}^nF(h_i)^{1-\ell}\right]f_n^{[\ell]}(\chi)+\left[\sum\limits_{i=1}^{n+1}F(h_i)^{1-\ell}\right]K_{n+1}^{[\ell]}\left(\|\chi-{\cal X}_{n+1}\| \right)},$$
with

\begin{eqnarray}
\varphi_n^{[\ell]}({\chi})\label{phin}=\frac{\sum\limits_{i=1}^n\frac{Y_i}{F(h_i)^{\ell}}K\left( \frac{\|{\chi}-{\cal X}_i\|}{h_i}\right)}{\sum\limits_{i=1}^nF(h_i)^{1-\ell}}, ~
f_n^{[\ell]}({\chi})\label{fn}=\frac{\sum\limits_{i=1}^n\frac{1}{F(h_i)^{\ell}}K\left( \frac{\|{\chi}-{\cal X}_i\|}{h_i}\right)}{\sum\limits_{i=1}^nF(h_i)^{1-\ell}},
\end{eqnarray}
and $K_i^{[\ell]}(\cdot):=\frac{1}{F(h_i)^{\ell}\sum\limits_{j=1}^iF(h_j)^{1-\ell}}K\left(\frac{\cdot}{h_i}\right)$.
%More precisely, $r_n^{[\ell]}({\chi})$ is  the adaption to the  functional model of the  finite-dimensional recursive family of estimators introduced by \cite{Amiri1}, which includes the famous ones, say recursive ($\ell=0$) and semi recursive ($\ell=1$) estimators.
The recursive property of this class of regression estimators offers many  advantages
and is clearly useful in sequential investigations and also for a large sample size.  Indeed,
%addition of a new observation means the non-recursive estimators must be entirely recomputed.  This
this kind of estimators are of  easy implementation and interpretation,  fast to compute and they do not require extensive
storage of data.
%Besides, we are required to store extensive data in order to calculate them.
The weak and strong consistency of this family of estimators was  studied by Amiri et al. \cite{Amiri2} in the framework of the independent case.

\section{Assumptions and main results}\label{resultat}
\subsection{Assumptions}
In the same spirit as Masry \cite{Masry2005}, we suppose throughout the paper   the existence of  nonnegative functions  $f_1$ and $\phi$ such that $\phi(0)=0$ and $F(h)=\Pb\left[\|\chi-\X\|\leq h\right]=\phi(h)f_1(\chi),$  for $h$ on a neighborhood  of zero. Then  $\phi$ is  an  increasing
function of $h$ and   $\phi(h)\to 0$ as $h\to 0$. The function  $f_1$  is referred to as  a  functional probability density (see Gasser et al. \cite{Gasser1998} for more details). We will assume that the following  assumptions hold.
\begin{enumerate}

\item [(H1)] The operators  $ r $ and $ \sigma_\varepsilon^2$ are continuous on a neighborhood of ${\chi}$.  Moreover, the function \\  $\zeta(t):=\Eb\left[\{r({\cal X})-r(\chi)\}~/~\|{\cal X}-\chi\|=t\right]$ is assumed to be derivable at $t=0$.
\item [(H2)] $K$ is nonnegative bounded kernel with support on the compact $[0,1]$ such that $\inf\limits_{t\in [0,1]}K(t)>0$.
\item [(H3)] For any $s\in[0,1], \tau_h(s):=\frac{\phi(hs)}{\phi(h)}\rightarrow\tau_0(s)<\infty$ as $h\rightarrow0$.
\item [(H4)]
\begin{enumerate}
 \item[(i)]$h_n\downarrow 0, ~n\phi(h_n)\rightarrow\infty$,  $A_{n,\ell}:=\displaystyle\frac{1}{n}\sum\limits_{i=1}^n\frac{h_i}{h_n}\left[\frac{\phi(h_i)}{\phi(h_n)}\right]^{1-\ell}\rightarrow\alpha_{[\ell]}>0$ as $n \to \infty$.
\item[(ii)] $\forall r\leq2$, $B_{n,r}:=\displaystyle\frac{1}{n}\sum\limits_{i=1}^n\left[\frac{\phi(h_i)}{\phi(h_n)}\right]^r\rightarrow\beta_{[r]}>0,$ as $n \to \infty$.
        \item[(iii)]  For any $\mu>0$, $\lim\limits_{n\to \infty}\dfrac{(\ln n)^{3+\frac{2}{\mu}}}{n\phi(h_n)}=0$.
\end{enumerate}
 \item [(H5)]
\begin{enumerate}
\item [(i)] $(\X_t)_{t\in \mathbb{N}}$ is a strong mixing process with $\alpha_{\X}(k)\leq ck^{-\rho}$, $k\geq 1$, for some $c>0$ and $\rho>2$.
%\item [(ii)] $\Pb\left(\X \in \mathcal{B}(\chi,h)\right)=F(h)=c\phi(h)f_1(\chi)$ where $\phi(h)\to 0$ as $h\to 0$ and $f_1$ is a nonnegative function.
\item[(ii)] There exist  non negative functions  $\psi$ and $f_2$  such that $\psi(h)\to 0$ as $h \to  0$, the ratio $\dfrac{\psi(h)}{\phi(h)^2}$ is bounded and
%  such that  %for any $i\neq j$,
 $\sup\limits_{i\neq j}\Pb\left[(\X_i,\X_j)\in \mathcal{B}(\chi,h_i)\times \mathcal{B}(\chi,h_j)\right]\leq \psi(h_i)\psi(h_j)f_2(\chi).$ %where $f_2$ is a nonnegative function.
 \end{enumerate}
 \item [(H6)]There exist   $\lambda>0$ and $\mu>0$  such that $\Eb\left[\exp\left(\lambda |Y|^\mu\right)\right]<\infty.$
 \end{enumerate}

 Since this paper is a generalization to dependent case  of the results  in Amiri et al.  \cite{Amiri2}, several of the assumptions are the same as those used in the earlier reference. The reader is then referred to this last for more comments on assumptions. Let us mention that the decrease of the sequence $(h_n)_{n\in \mathbb{N}}$  is   particular to the recursive estimators and for dependent data. The technical condition (H4)(iii) is unrestrictive and is easily  satisfied by the popular choices of $\phi$ and $h_n$ given by $\phi(h_n)\sim n^{-\xi},$ with $0<\xi<1$. Assumption (H5)(i) is  the classically strong   mixing condition which,  is well known to be  satisfied  by linear or stationary ARMA processes.  In order to simplify the presentation, we assume the strong  mixing coefficient  to be arithmetic, but the main results can be obtained under  several  conditions on this coefficient. Assumption (H5)(ii) plays a crucial role in our calculus, when we show the negligibility of some covariance terms. It has been used by Masry  \cite{Masry2005} in the   non recursive  case. Finally,  as developped  in Amiri et al.  \cite{Amiri2}, assumption (H6) implies that
\begin{eqnarray}\label{moment-expo}\Eb\left( \max\limits_{1\leq i\leq n}|Y_i|^p\right) = O[(\ln n)^{p/\mu}], \forall  p \geq1,n \geq 2.\end{eqnarray}
\subsection{Main results }
For convenience, let us introduce the following notations:
\begin{eqnarray*}
M_0&=&K(1)-\int_{0}^1(sK(s))'\tau_0(s)ds,~~
M_1=K(1)-\int_{0}^1K'(s)\tau_0(s)ds\\
M_2&=&K^2(1)-\int_{0}^1(K^2(s))'\tau_0(s)ds.\\
\end{eqnarray*}
In the following theorem, we establish the almost sure convergence of  the proposed recursive kernel estimator of the regression function.
\begin{theorem}\label{Cv_ps_rnl}
     Assume that (H1)-(H6) hold.
 If   $ \displaystyle\lim\limits_{n\rightarrow+\infty}nh_n^{2}=0$, then
$$\limsup_{n\rightarrow\infty} \left[\frac{n\phi(h_n)}{\ln~n}\right]^{1/2}\left[r_n^{[\ell]}(\chi)-r(\chi)\right]\leq \frac{2}{M_1}\left[ 1+V_\ell(\chi) \right]~ a.s.$$
where
\begin{eqnarray}\label{V}
V_\ell(\chi)= \frac{\beta_{[1-2\ell]}}{\beta_{[1-\ell]}^2}\frac{\sigma_\varepsilon^2(\chi)}{f_1(\chi)}M_2,
\end{eqnarray}
for all $\chi$ such that $f_1(\chi)>0$.
\end{theorem}

Theorem \ref{Cv_ps_rnl} is an extension of  Ferraty and Vieu's  \cite{Ferraty2004} result  on functional kernel-type  estimate to the general family of recursive estimators $r^{[\ell]}_n (\chi).$ A similar  result is also obtained by Ling and Wu \cite{Ling2012} for a truncated version  of the  Nadaraya-Watson type estimator, under the condition  $\Eb|Y |<\infty$, which is weaker than  assumption (H6). However,  Theorem \ref{Cv_ps_rnl}  establishes the rate of convergence with exact appear bound, while   Ling and Wu's \cite{Ling2012} result tells only  the rate of convergence in function of the variances of the  numerator and denominator of the estimator. As we will see in the proofs below, assumption (H6)  will be necessary, for the  study of the covariance terms and also  when we shall  prove the cancellation of the residual term between the estimator and its  truncated version.   Finally, let us mention that
compared with the result in Amiri et al.  \cite{Amiri2}, as in the multivariate framework, it is difficult to obtain the optimal rate  $ \left[\frac{n\phi(h_n)}{\ln\ln~n}\right]^{1/2}$ in the dependent case. \\

\noindent The mean square  error  of $r_n^{[\ell]}(x)$ is given in Theorem \ref{mse} below.

\begin{theorem}\label{mse}
Under assumptions (H1)-(H6),
\begin{eqnarray*}
\Eb\left[\left(r_n^{[\ell]}({\chi})-r(\chi)\right)^2\right]=\left[\left(\frac{\zeta'(0)\alpha_{[\ell]}M_0h_n}{\beta_{[1-\ell]}M_1}\right)^2+\frac{\beta_{[1-2\ell]}M_2\sigma_\varepsilon^2(\chi)}{\beta_{[1-\ell]}^2M_1^2f_1(\chi)n\phi(h_n)}\right]\left[1+o(1)\right]
\end{eqnarray*}
for all $\chi$ such that $f_1(\chi)>0$.\end{theorem}
Theorem \ref{mse} is an extension to functional data of the result of Amiri \cite{Amiri1} in finite dimensional setting. Also, our result generalizes the works of Bosq and Cheze-Payaud \cite{bosq-cheze1999} to functional and recursive setting. Finally, in counterpart of the almost sure convergence, Theorem \ref{mse} gives the same rate of convergence and asymptotic  constants as those obtained for  the iid case in  Amiri et al. \cite{Amiri2}.
\section{Proofs}\label{preuve}
In the sequel,  and through  the paper, $c$ will denote a constant whose value is unimportant and may vary from line to line. Also,
we  set \begin{eqnarray*}
K_i(\chi)=K\left(\frac{\|\chi-\X_i\|}{h_i}\right).
\end{eqnarray*}
Finally,  for convenience we will use   the following decomposition
\begin{eqnarray}\label{egalite1}
 r_n^{[\ell]}(\chi)-r(\chi)=\frac{\tilde{\varphi}_n^{[\ell]}(\chi)-r(\chi)f_n^{[\ell]}(\chi)}{f_n^{[\ell]}(\chi)}+\frac{\varphi_n^{[\ell]}(\chi)-
\tilde{\varphi}_n^{[\ell]}(\chi)}{f_n^{[\ell]}(\chi)},
\end{eqnarray}
%o: $\varphi_n^{[\ell]}({\chi})=\frac{1}{\sum\limits_{i=1}^nF(h_i)^{1-\ell}}\sum\limits_{i=1}^n\dfrac{Y_i}{F(h_i)^{\ell}}K\left( \dfrac{\|{\chi}-{\cal X}_i\|}{h_i}\right),~f_n^{[\ell]}({\chi})=\frac{1}{\sum\limits_{i=1}^nF(h_i)^{1-\ell}}\sum\limits_{i=1}^n\dfrac{1}{F(h_i)^{\ell}}K\left( \dfrac{\|{\chi}-{\cal X}_i\|}{h_i}\right)$ et
where $\tilde{\varphi}_n^{[\ell]}(\chi)$ is a truncated version of  $\varphi_n^{[\ell]}(\chi)$ defined by
 \begin{equation}\label{phi_nl_tilde}\tilde{\varphi}_n^{[\ell]}({\chi})=\frac{1}{\sum\limits_{i=1}^nF(h_i)^{1-\ell}}\sum\limits_{i=1}^n\dfrac{Y_i}{F(h_i)^{\ell}}\mathds{1}_{\left\{\left|Y_i\right|\leq b_n\right\}}K\left( \dfrac{\|{\chi}-{\cal X}_i\|}{h_i}\right),\end{equation}
$b_n$ being a sequence of real numbers which goes to $+\infty$ as $n\rightarrow\infty.$
\subsection{Preliminary lemmas}
In order to prove the main  results, we need the following lemmas.

\begin{lemma}\label{Biais1}
Under assumptions (H1)-(H4), we have
\begin{eqnarray*}
\frac{\Eb\left[\varphi_n^{[\ell]}({\chi})\right]}{\Eb\left[f_n^{[\ell]}({\chi})\right]}-r(\chi)=h_n\zeta'(0)\frac{\alpha_{[\ell]}}{\beta_{[1-\ell]}}\frac{M_0}{M_1}\left[1+o(1)\right].
\end{eqnarray*}
\end{lemma}
{\it Proof.} See   Amiri et al. \cite{Amiri2}, since the bias term is not depending to the mixing structure.
%Lemma \ref{lem3} is proved in \cite{Amiri2} (See the proof of Lemma 6)..
 \hfill $\square$
\begin{lemma}\label{Var}
Under assumptions (H1)-(H6), we have
\begin{eqnarray*}
\text{Var}\left[f_n^{[\ell]}({\chi})\right]&=&\frac{\beta_{[1-2\ell]}}{\beta_{[1-\ell]}^2}\frac{M_2}{f_1(\chi)}\frac{1}{n\phi(h_n)}\left[1+o(1)\right];\\
\text{Var}\left[\varphi_n^{[\ell]}({\chi})\right]&=&\frac{\beta_{[1-2\ell]}}{\beta_{[1-\ell]}^2}\left[r^2(\chi)+\sigma_\epsilon^2(\chi)\right]\frac{M_2}{f_1(\chi)}\frac{\left[1+o(1)\right]}{n\phi(h_n)};\\
\text{Cov}\left[f_n^{[\ell]}({\chi}),\varphi_n^{[\ell]}({\chi})\right]&=&\frac{\beta_{[1-2\ell]}}{\beta_{[1-\ell]}^2}r(\chi)\frac{M_2}{f_1(\chi)}\frac{\left[1+o(1)\right]}{n\phi(h_n)},
\end{eqnarray*}
for all  $\chi$ such that $f_1(\chi)>0.$
\end{lemma}
{\it Proof.}
%\subsection{Proof of Lemma \ref{Var}}
The variance term of $f_n^{[\ell]}(\chi)$ can be decomposed in variance and covariance terms as
\begin{eqnarray}\label{varfnl}
\text{Var}(f_n^{[\ell]}(\chi))&=& \left[\sum_{i=1}^nF(h_i)^{1-\ell}\right]^{-2}\left(\sum_{i=1}^nA_{i,i}+\sum_{i\neq j}A_{i,j}\right) :=  F_1+F_2,
\end{eqnarray}
where  for any integers $i$ and $j$,
%\begin{eqnarray*}
$A_{i,j}=F(h_i)^{-\ell}F(h_j)^{-\ell}\text{Cov}\left(K_i(\chi),K_j(\chi)\right).$
%\end{eqnarray*}
Noting that the principal term $F_1$ in the right-hand side of \eqref{varfnl} corresponds to the variance term of $f_n^{[\ell]}$ in the independent case (see Amiri et al. \cite{Amiri2} fore more details), and is given by
\begin{eqnarray*}
n\phi(h_n)F_1=\frac{\beta_{[1-2\ell]}}{\beta_{[1-\ell]}^2}\frac{M_2}{f_1(\chi)}\left[1+o(1)\right].
\end{eqnarray*}
Now, let us establish that the covariance term $F_2$ is negligible. To this end, let $c_n$ be a sequence of real numbers tending to $\infty$ as $n \to \infty$. We can write
\begin{eqnarray}\label{F}
F_2&\leq &\frac{2\left(\sum_{k=1}^{c_n}\sum_{p=1}^n|A_{k+p,p}|+\sum_{k=c_n+1}^{n-1}\sum_{p=1}^n|A_{k+p,p}|\right)}{\left[\sum_{i=1}^nF(h_i)^{1-\ell}\right]^{2}} :=F_{21}+F_{22}.
\end{eqnarray}
From assumptions (H2) and (H5)(ii),  we have for any $i\neq j$
\begin{eqnarray}\label{KiKj}
\Eb\left[K_i(\chi)K_j(\chi)\right]& =&\int_{[0,1]\times[0,1]} K(u)K(v)d\Pb^{\left(\frac{\|\chi-\X_i\|}{h_i},\frac{\|\chi-\X_j\|}{h_j}\right)}(u,v)\nonumber\\
& \leq & \|K\|_\infty^2\Pb\left[\|\chi-\X_i\|\leq h_i,\|\chi-\X_j\|\leq h_j\right]\nonumber\\
& \leq & c\psi(h_i)\psi(h_j).
\end{eqnarray}
Note that, from the proof of Lemma 2  in  Amiri et al. \cite{Amiri2} we can write
\begin{eqnarray*}
\Eb\left[K_i(\chi)\right]= \phi(h_i)f_1(\chi)\left[K(1)-\int_0^1K'(s)\tau_{h_i}(s)ds\right],
\end{eqnarray*}
so, we get
\begin{eqnarray}\label{inegal_cov}
\left|\cov(K_i(\chi),K_j(\chi))\right| &\leq &c\left[\psi(h_i)\psi(h_j)+\phi(h_i)\phi(h_j)\right].
\end{eqnarray}
Hence, we deduce that
\begin{eqnarray*}
F_{21}& \leq& \frac{c \sum_{k=1}^{c_n}\sum_{p=1}^n\left[\frac{\psi(h_{k+p})\psi(h_p) }{\phi(h_{k+p})^{\ell}\phi(h_p)^{\ell}}+\phi(h_{k+p})^{1-\ell}\phi(h_p)^{1-\ell}\right]}{\left[\sum_{i=1}^n\phi(h_i)^{1-\ell}\right]^{2}}: =F_{211}+F_{212}.
\end{eqnarray*}
Now, Assumption (H5)(ii)  ensures that the ratio $\psi(h_i)/\phi(h_i)$ is bounded  and since $\phi$ is increasing, we get
\begin{eqnarray*}
F_{211} & \leq& c\left[\sum_{i=1}^n\phi(h_i)^{1-\ell}\right]^{-2}\sum_{k=1}^{c_n}\sum_{p=1}^n\phi(h_p)^{2-2\ell}\leq
c\frac{B_{n,2-2\ell}}{B_{n,1-\ell}^2}\frac{c_n}{n}.
\end{eqnarray*}
Hence,
\begin{eqnarray}\label{F211}
n\phi(h_n)F_{211}=O\left(\phi(h_n)c_n\right).
\end{eqnarray}
Now for the second term $F_{212}$, again, using the fact that $\phi$ is an  increasing  function, we get
\begin{eqnarray*}
F_{212} & \leq & c\left[\sum_{i=1}^n\phi(h_i)^{1-\ell}\right]^{-2}\sum_{k=1}^{c_n}\sum_{p=1}^n\phi(h_p)^{2-2\ell} \leq c\frac{B_{n,2-2\ell}}{B_{n,1-\ell}^2}\frac{c_n}{n},
\end{eqnarray*}
so that
\begin{eqnarray}\label{F212}
n\phi(h_n)F_{212}=O\left(\phi(h_n)c_n\right).
\end{eqnarray}
From \eqref{F211} and \eqref{F212} we deduce
\begin{eqnarray}\label{F21}
n\phi(h_n)F_{21}=O\left(\phi(h_n)c_n\right).
\end{eqnarray}
Next, for the second term $F_{22}$ in \eqref{F}, we have from Billingsley's inequality,
\begin{eqnarray*}
F_{22} & \leq & c\left[\sum_{i=1}^n\phi(h_i)^{1-\ell}\right]^{-2}\sum_{k=c_n+1}^{n-1}\sum_{p=1}^nk^{-\rho}\phi(h_{k+p})^{-\ell}\phi(h_p)^{-\ell}\\& \leq &
c\left[\sum_{i=1}^n\phi(h_i)^{1-\ell}\right]^{-2}\frac{c_n^{1-\rho}}{\rho-1}\sum_{p=1}^n\phi(h_n)^{-\ell}\phi(h_p)^{-\ell}\leq
c\frac{B_{n,-\ell}}{B_{n,1-\ell}^2}\frac{c_n^{1-\rho}}{n\phi(h_n)^2}.
\end{eqnarray*}
Therefore
\begin{eqnarray}\label{F22}
n\phi(h_n)F_{22}=O\left(\frac{c_n^{1-\rho}}{\phi(h_n)}\right).
\end{eqnarray}
If we choose  $c_n=\lfloor \phi(h_n)^{-\frac{2}{\rho}}\rfloor$, we deduce from \eqref{F21} and \eqref{F22} that
\begin{eqnarray*}
n\phi(h_n)F_2&=&O\left(\phi(h_n)^{\frac{\rho-2}{\rho}}\right)=o(1) \text{  as long  as }\rho>2,
\end{eqnarray*}
and the first part of Lemma \ref{Var} follows. Now, as in the proof of the first part of Lemma \ref{Var}, the variance term of $\varphi_n^{[\ell]}$ is decomposed as follows
\begin{eqnarray*}
\text{Var}(\varphi_n^{[\ell]}(\chi))=\left[\sum_{i=1}^nF(h_i)^{1-\ell}\right]^{-2}\left(\sum_{i=1}^nA_{i,i}+\sum_{i\neq j}A_{i,j}\right)  :=  I_1+I_2,
\end{eqnarray*}
where here $A_{i,j}$ denotes for any integers $i$ and $j$ as follows
%\begin{eqnarray*}
$A_{i,j}=F(h_i)^{-\ell}F(h_j)^{-\ell}\text{Cov}\left(Y_iK_i(\chi),Y_jK_j(\chi)\right).$
%\end{eqnarray*}
The study of the term $I_1$ is treated in the same manner as in the independent case (see Amiri et al.  \cite{Amiri2} for more details) which gives
\begin{eqnarray*}
n\phi(h_n)I_1=\frac{\beta_{[1-2\ell]}}{\beta_{[1-\ell]}^2}\frac{r(\chi)}{f_1(\chi)}M_2\left[1+o(1)\right].
\end{eqnarray*}
Now, for the second term $I_2$, we always consider a sequence of real numbers $c_n$ which goes to $\infty$ as $n\to\infty$ and we write
\begin{eqnarray}\label{I}
I_2&\leq& \frac{2\left(\sum_{k=1}^{c_n}\sum_{p=1}^n|A_{k+p,p}|+\sum_{k=c_n+1}^{n-1}\sum_{p=1}^n|A_{k+p,p}|\right)}{\left[\sum_{i=1}^nF(h_i)^{1-\ell}\right]^{2}}:=I_{21}+I_{22}.
\end{eqnarray}
The term $I_{22}$ is treated exactly as $F_{22}$ in the proof of the first part of this Lemma previously, by substituting the Billingsley lemma with the Davydov lemma. Then, setting $b_n=(\delta\ln n)^{1/\mu}$, using \eqref{moment-expo} and with the help of (H5), we get
\begin{eqnarray*}
I_{22}&\leq& c(\ln n)^{2/\mu}\left[\sum_{i=1}^n\phi(h_i)^{1-\ell}\right]^{-2}\frac{c_n^{1-\rho/2}}{(\rho/2)-1}\phi(h_n)^{-2\ell}\sum_{p=1}^n\left(\frac{\phi(h_p)}{\phi(h_n)}\right)^{-\ell} \\&\leq &
c (\ln n)^{2/\mu}\frac{B_{n,-\ell}}{B_{n,1-\ell}^2}\frac{1}{n\phi(h_n)^2}c_n^{1-\rho/2}.
\end{eqnarray*}
Therefore,
\begin{eqnarray}\label{I22}
n\phi(h_n)I_{22}&=&O\left((\ln n)^{2/\mu}c_n^{1-\rho/2}(\phi(h_n))^{-1}\right).
\end{eqnarray}
For the second term $I_{21}$, observe that for any integers $i$ and $j$,
\begin{eqnarray*}
\left|\cov(Y_iK_i(\chi),Y_jK_j(\chi))\right| & \leq& \left|\Eb\left[Y_iY_jK_i(\chi)K_j(\chi)\right]\right|+\left|\Eb\left[Y_iK_i(\chi)\right]\right|\left|\Eb\left[Y_jK_j(\chi)\right]\right|.
\end{eqnarray*}
Now, from assumptions (H1) and (H2) and conditioning on $\X$, one have
\begin{eqnarray*}
\Eb\left[Y_iK_i(\chi)\right]&=&
M_1F(h_i)\left[r(\chi)+\gamma_i\right]\leq  c\phi(h_i),
\end{eqnarray*}
where $\gamma_i$ goes to zero as $i\to\infty$.
Using Cauchy-Schwartz' inequality, choosing $b_n=(\ln n)^{1/\mu}$ and \eqref{moment-expo}, we get
\begin{eqnarray*}
\left|\cov(Y_iK_i(\chi),Y_jK_j(\chi))\right| & \leq& \Eb^{1/2}\left[Y_i^2Y_j^2\right]\Eb^{1/2}\left[K_i^2(\chi)K_j^2(\chi)\right]\\&&+\left|\Eb\left[Y_iK_i(\chi)\right]\right|\left|\Eb\left[Y_jK_j(\chi)\right]\right|\\ &\leq & c\left[(\ln n)^{2/\mu}\psi(h_i)^{1/2}\psi(h_j)^{1/2}+\phi(h_i)\phi(h_j)\right].
\end{eqnarray*}
The rest of the proof for $I_{21}$ is the same as the one for $F_{21}$ which implies that
%\begin{eqnarray*}
$I_{21}\leq  c\frac{c_n}{n}\left[(\ln n)^{2/\mu}+1\right].$
%\end{eqnarray*}
Hence,
\begin{eqnarray}\label{I21}
n\phi(h_n)I_{21}=O\left(c_n\phi(h_n)(\ln n)^{2/\mu}\right),
\end{eqnarray}
and the result of the second part of Lemma \ref{Var} follows from \eqref{I22} and \eqref{I21}  with the choice $c_n=\lfloor \phi_n^{-\frac{4}{\rho}}\rfloor$. Next, to treat the last part of Lemma \ref{Var}, it suffices to decompose the  term  $n\phi(h_n)\cov\left(\phi_n^{[\ell]}(\chi),f_n^{[\ell]}(\chi)\right)$ by the principal and covariance terms and use the same procedure as in the proof of the second part of Lemma \ref{Var}. \hfill $\square$

\begin{lemma} \label{lem1} Set $$N=\left[\frac{n\phi(h_n)}{\ln~ n}\right]^{1/2}
\left\{\tilde{\varphi}_n^{[\ell]}(\chi)-r(\chi)f_n^{[\ell]}(\chi)-\Eb\left[\tilde{
\varphi}_n^{[\ell]}(\chi)-r(\chi)f_n^{[\ell]}(\chi)\right]\right\}, $$
where $\tilde{\varphi}$ is defined in \eqref{phi_nl_tilde}. Under assumptions (H1)-(H6), we have
\begin{eqnarray*}
\label{I1reg}\varlimsup_{n\rightarrow\infty}N\leq 2\left[1+V_{\ell}(\chi)\right] \text{ a.s.},
\end{eqnarray*}
where $V_\ell$ is defined in \eqref{V}.
\end{lemma}
%For the first term of the right hand side of \eqref{decomposition11},  we have the following lemma.
{\it Proof.}
%\subsubsection{Proof of Lemma \ref{lem1}}
Set
%\begin{eqnarray*}
$W_{n,i}=\dfrac{K_i(\chi)\left[Y_i\mathds{1}_{\left\{|Y_i|\leq b_n\right\}}-r(\chi)\right]}{f_1(\chi)B_{n,1-\ell}\phi(h_n)^{1-\ell}\phi(h_i)^\ell},\text{  where }Z_{n,i}= W_{n,i}-\text{E}W_{n,i}.$
%\end{eqnarray*}
To prove Lemma \ref{lem1}, we use the blocks decomposition technique.
Let $p_n$ and $q_n$ be some sequences of real numbers defined by
$p_n=\lfloor p_0\ln n\rfloor$ with $p_0>0$ and $q_n=\lfloor \frac{n}{2p_n}\rfloor$. Set
\begin{eqnarray*}
S_n'=\sum_{j=1}^{q_n}V_n(2j-1), \ \ S_n''=\sum_{j=1}^{q_n}V_n(2j) \text{ and } S_n'''=\frac{1}{n}\sum_{k=2p_nq_n+1}^nZ_{n,k},
\end{eqnarray*}
with
%\begin{eqnarray*}
$V_n(j)=\frac{1}{n}\sum\limits_{k=(j-1)p_n+1}^{jp_n}Z_{n,k}, \ \ j=1,\ldots,2q_n.$
%\end{eqnarray*}
Then we have
%\begin{eqnarray*}
$N=S_n'+S_n''+S_n'''.$
%\end{eqnarray*}
Observe that the third term $S_n'''$ is negligible so that, to prove the strong consistency of $N$, it suffices to check the almost sure convergence for $S_n'+S_n''$. For any $\varepsilon>0$,
\begin{eqnarray*}
\Pb\left(|S_n'+S_n''|>\varepsilon\right)\leq \Pb\left(|S_n'|>\frac{\varepsilon}{2}\right)+\Pb\left(|S_n''|>\frac{\varepsilon}{2}\right).
\end{eqnarray*}
We just treat $S_n'$, the term $S_n''$ being similar. Since $K$ is bounded and $\phi$ is non decreasing, we get for $n$ large enough
%\begin{eqnarray*}
$\left|V_n(j)\right| \leq \frac{2\|K\|_{\infty}p_nb_n}{f_1(\chi)B_{n,1-\ell}n\phi(h_n)}.$
%\end{eqnarray*}
Using Rio's \cite{Rio2000} coupling lemma, the random variables $V_n(j)$ can be approximated by independent and identically distributed random variables $V_n^*(j)$ such that
%\begin{eqnarray*}
$$\Eb\left|V_n(2j-1)-V_n^*(2j-1)\right| \leq \frac{4\|K\|_\infty p_nb_n}{f_1(\chi)B_{n,1-\ell}n\phi(h_n)}\alpha(p_n).$$
%\end{eqnarray*}
Since $p_nq_n\leq n$, it follows that
\begin{eqnarray*}
\sum_{j=1}^{q_n}\Eb\left|V_n(2j-1)-V_n^*(2j-1)\right| &\leq &\frac{4\|K\|_\infty p_nq_n b_n}{f_1(\chi)B_{n,1-\ell}n\phi(h_n)}\alpha(p_n)\\ &\leq &\frac{4\|K\|_\infty b_n}{f_1(\chi)B_{n,1-\ell}\phi(h_n)}\alpha(p_n).
\end{eqnarray*}
Therefore, for any $\varepsilon, \kappa>0$, Markov's inequality leads to
\begin{eqnarray}\label{Inegalite1}
\Pb\left(\left|\sum_{j=1}^{q_n}\left[V_n(2j-1)-V_n^*(2j-1)\right]\right|>\frac{\varepsilon\kappa}{2(1+\kappa)}\right)&\leq& \frac{8(1+\kappa)}{\varepsilon\kappa}\frac{\|K\|_\infty b_n\alpha(p_n)}{f_1(\chi)B_{n,1-\ell}\phi(h_n)}\nonumber\\ & \leq &
\frac{8(1+\kappa)}{\varepsilon\kappa}\frac{\|K\|_\infty b_n\gamma e^{-\rho p_0\ln n}}{f_1(\chi)B_{n,1-\ell}\phi(h_n)}.
\end{eqnarray}
Next setting $\varepsilon_n=\varepsilon\sqrt{\dfrac{\ln n}{n\phi(h_n)}}$ with $\varepsilon>0$ and $\lambda_n=\sqrt{n\phi(h_n)\ln n}$, we have from (H4)(iii),
%\begin{eqnarray*}
$\left|\lambda_nV_n^*(j)\right|  \leq   \frac{2\|K\|_\infty p_nb_n}{f_1(\chi)B_{n,1-\ell}}\sqrt{\dfrac{\ln n}{n\phi(h_n)}} \rightarrow 0,$
%\end{eqnarray*}
therefore for $n$ large enough,
%\begin{eqnarray*}
$\left|\lambda_nV_n^*(j)\right|\leq \frac{1}{2}.$
%\end{eqnarray*}
It follows that
\begin{eqnarray*}
\exp\left\{\pm \lambda_nV_n^*(j)\right\} \leq 1\pm \lambda_nV_n^*(j)+\left[\lambda_nV_n^*(j)\right]^2.
\end{eqnarray*}
From Markov's inequality, we get
\begin{eqnarray*}
\Pb\left(\left|\sum_{j=1}^{q_n}V_n^*(2j-1)\right|>\frac{\varepsilon_n}{2(1+\kappa)}\right)&\leq &\Pb\left(\exp\left[\sum_{j=1}^{q_n}\lambda_nV_n^*(2j-1)\right]>\exp\left[\frac{\lambda_n\varepsilon_n}{2(1+\kappa)}\right]\right) \\ & &
+ \Pb\left(\exp\left[-\sum_{j=1}^{q_n}\lambda_nV_n^*(2j-1)\right]>\exp\left[\frac{\lambda_n\varepsilon_n}{2(1+\kappa)}\right]\right)\\
& \leq &
2 \exp\left[\frac{-\lambda_n\varepsilon_n}{2(1+\kappa)}+\lambda_n^2\sum_{j=1}^{q_n}\Eb V_n^{*2}(2j-1)\right].
\end{eqnarray*}
Since,
%\begin{eqnarray*}
$\sum\limits_{j=1}^{q_n}\Eb V_n^{*2}(2j-1)\leq  \dfrac{1}{n^2}\left[\sum\limits_{k=1}^n\text{Var}(Z_{n,k})+\sum\limits_{k\neq k'}\text{Cov}(Z_{n,k},Z_{n,k'})\right],$
%\end{eqnarray*}
we will assume for the moment that

\begin{eqnarray}\label{moment1}
\dfrac{\phi(h_n)}{n}\sum_{k=1}^n\text{Var}(Z_{n,k}) &= &V_\ell(\chi)\left[1+o(1)\right] \\ \dfrac{\phi(h_n)}{n}\sum_{k\neq k'}\text{Cov}(Z_{n,k},Z_{n,k'}) &= &o(1),\label{moment2}
\end{eqnarray}
where $V_\ell$ is defined in \eqref{V}.
It follows from \eqref{moment1} and \eqref{moment2} that, for $n$ large enough,
%\begin{eqnarray*}
$$\lambda_n^2\sum_{j=1}^{q_n}\Eb V_n^{*2}(2j-1)\leq V_\ell(\chi)\ln n\left[1+o(1)\right].$$
%\end{eqnarray*}
Therefore
\begin{eqnarray}\label{Inegalite2}
\Pb\left(\left|\sum_{j=1}^{q_n}V_n^*(2j-1)\right|>\frac{\varepsilon_n}{2(1+\kappa)}\right)\leq 2 e^{\left[\frac{-\varepsilon}{2(1+\kappa)}+V_\ell(\chi)(1+o(1))\right]\ln n}
\end{eqnarray}
Now, combining \eqref{Inegalite1} and \eqref{Inegalite2}, we get
\begin{eqnarray}\label{sn'}
\Pb\left(|S_n'|>\frac{\varepsilon_n}{2}\right)&\leq& \Pb\left(\left|\sum_{j=1}^{q_n}V_n(2j-1)-V_n^*(2j-1)\right|>\frac{\varepsilon_n\kappa}{2(1+\kappa)}\right)\nonumber\\ & & +\Pb\left(\left|\sum_{j=1}^{q_n}V_n^*(2j-1)\right|>\frac{\varepsilon_n\kappa}{2(1+\kappa)}\right)\nonumber\\ & \leq & \gamma \frac{8(1+\kappa)}{\varepsilon\kappa}\frac{\|K\|_\infty b_nn^{1-p_0\rho}}{f_1(\chi)B_{n,1-\ell}\sqrt{\ln n}\sqrt{n\phi(h_n)}}\nonumber\\ & & +
2 \exp\left\{\left[-\frac{\varepsilon}{2(1+\kappa)}+V_\ell(\chi)\right]\ln n\right\}.
\end{eqnarray}
Next, with the choice of $b_n=(\delta \ln n)^{1/\mu}$, the conclusion follows from the application of the Borel-Cantelli's lemma whenever $p_0>\frac{2}{\rho}$ and $\varepsilon>2(1+\kappa)\left[1+V_\ell(\chi)\right]$, which implies that
\begin{eqnarray*}
\varlimsup_{n\rightarrow\infty}\left[\frac{n\phi(h_n)}{\ln~ n}\right]^{1/2}
N\leq 2(1+\kappa)\left[ 1+V_{\ell}(\chi)\right] \text{ a.s.},
\end{eqnarray*}
for all positive $\kappa$ and Lemma \ref{lem1} follows.
To complete the proof, let us prove \eqref{moment1} and \eqref{moment2}.  We can write
\begin{eqnarray*}
\dfrac{\phi(h_n)}{n}\sum_{k=1}^n\text{Var}(Z_{n,k})=\dfrac{\sum_{k=1}^n\phi(h_k)^{-2\ell}\text{Var}\left(K_i(\chi)\left[Y_i\mathds{1}_{\left\{|Y_i|\leq b_n\right\}}-r(\chi)\right] \right)}{f_1^2(\chi)B_{n,1-\ell}^2n\phi(h_n)^{1-2\ell}}.
\end{eqnarray*}
Following the same lines of the proof of Lemma 5 in Amiri et al.  \cite{Amiri2}, one can prove that
\begin{eqnarray*}
\sum_{k=1}^n\phi(h_k)^{-2\ell}\text{Var}\left(K_k(\chi)\left[Y_k\mathds{1}_{\left\{|Y_k|\leq b_n\right\}}-r(\chi)\right] \right)\sim n\phi(h_n)^{1-2\ell}\beta_{[1-2\ell]}\sigma^2_\varepsilon(\chi)M_2,
\end{eqnarray*}
therefore \eqref{moment1} follows.
%
%ICI
%We can write
%\begin{eqnarray*}
%\dfrac{\phi(h_n)}{n}\sum_{k=1}^n\text{Var}(Z_{n,k})=\frac{V_n}{B_{n,1-\ell}^2 f_1(\chi)n\phi(h_n)^{1-2\ell}},
%\end{eqnarray*}
%where $V_n$ is defined in  \cite{Amiri2}. Therefore, \eqref{moment1} follows from the equivalence
%$$V_n\sim nF(h_n)^{1-2\ell}\beta_{[1-2\ell]}\sigma^2_\varepsilon(\chi)M_2,  \text{  as } n\rightarrow+\infty.$$ (we refer to Lemma 6  in  \cite{Amiri2} for details).
%ICI
Next, about  the covariance term in \eqref{moment2}, for any integers $i\neq j$, let
\begin{eqnarray*}
A_{i,j}= F(h_i)^{-\ell} F(h_j)^{-\ell}\left|\cov\left(K_i(\chi)\left[Y_i\mathds{1}_{\left\{|Y_i|\leq b_n\right\}}-r(\chi)\right],K_j(\chi)\left[Y_j\mathds{1}_{\left\{|Y_j|\leq b_n\right\}}-r(\chi)\right]\right)\right|.
\end{eqnarray*}
Then, we have

\begin{eqnarray*}
\frac{1}{n^2}\sum_{k\neq k'}\text{Cov}(Z_{n,k},Z_{n,k'}) &\leq& \dfrac{2\left[\sum_{k=1}^{c_n}\sum_{p=1}^n|A_{k+p,p}|+\sum_{k=c_n+1}^{n-1}\sum_{p=1}^n|A_{k+p,p}|\right]}{f_1^2(\chi)B_{n,1-\ell}^2n^2\phi(h_n)^{2-2\ell}}\\
&:=&J_1+J_2.
\end{eqnarray*}
Using Billingsley's inequality, one can prove that
%\begin{eqnarray*}
$n\phi(h_n)J_2=O\left(b_n^2c_n^{1-\rho}\phi(h_n)^{-1}\right).$
%\end{eqnarray*}
Next, since
\begin{eqnarray*}
|A_{k+p,p}|& \leq& \left(b_n+|r(\chi)|\right)^2\left[\Eb\left(K_{k+p}(\chi)K_p(\chi)\right)+\Eb(K_{k+p}(\chi))\Eb(K_{p}(\chi))\right]
\\ & \leq &  c\left(b_n+|r(\chi)|\right)^2\left[\psi(h_{k+p})\psi(h_p)+\phi(h_{k+p})\phi(h_p)\right]\phi(h_{k+p})^{-\ell}\phi(h_p)^{-\ell}.
\end{eqnarray*}
Therefore, as in the proof of the first part of Lemma \ref{Var}, we get
\begin{eqnarray*}
n\phi(h_n)J_1=O(b_n^2\phi(h_n)c_n),
\end{eqnarray*}
which together with the choice  $c_n=\lfloor \phi(h_n)^{-\frac{2}{\rho}}\rfloor$  imply  \eqref{moment2} as long as $\rho>2$.
 \hfill $\square$

\subsection{Proofs of the main results}

\subsubsection{Proof of Theorem \ref{Cv_ps_rnl}}

%Next, for any $\varepsilon>0$, we have for the residual term of (\ref{egalite1})
%$$P\left\{\left|\varphi_n^{[\ell]}(\chi)-\tilde{\varphi}_n^{[\ell]}(\chi)\right|>\varepsilon\left[\frac{\ln n}{ n\phi(h_n)}\right]^\frac{1}{2}\right\}
% \leq P\left(\bigcup\limits_{i=1}^n\left\{|Y_i|>b_n\right\}\right)
%  \leq\text{E}\left[e^{\lambda|Y|^\mu
% }\right]n^{1-\lambda\delta},$$
% where the last inequality follows by setting $b_n=(\delta\ln n)^\frac{1}{\mu},$  with the help of Markov's inequality.  Assumption ${\bf H5}$ ensures that for any $\varepsilon>0$, $$\sum_{n=1}^\infty P\left\{\left|\varphi_n^{[\ell]}(\chi)-\tilde{\varphi}_n^{[\ell]}(\chi)\right|>\varepsilon\left[\frac{ \ln n}{ n\phi(h_n)}\right]^\frac{1}{2}\right\}<\infty \text{ if  $\delta>\frac{2}{\lambda}$} ,$$   and the Borel-Cantelli's lemma implies that
Let us consider the decomposition (\ref{egalite1}). For the residual term $\frac{\varphi_n^{[\ell]}(\chi)-
\tilde{\varphi}_n^{[\ell]}(\chi)}{f_n^{[\ell]}(\chi)}$, following the same lines of proof in Amiri et al.  \cite{Amiri2} by replacing $\frac{\ln\ln n}{nF(h_n)}$ by $\frac{\ln n}{n\phi(h_n)}$, one can show that
 \begin{eqnarray}\label{conv_residual}
 \left[\frac{ n\phi(h_n)}{\ln n}\right]^{1/2}\left|\varphi_n^{[\ell]}(\chi)-\tilde{\varphi}_n^{[\ell]}(\chi)\right|\rightarrow 0 \text{ a.s,} \text{ when }n\rightarrow\infty.
 \end{eqnarray}

 For the principal term in \eqref{egalite1}, we can write
 \begin{eqnarray}
 \label{decomposition11}
 \tilde{\varphi}_n^{[\ell]}(\chi)-r(x)f_n^{[\ell]}(\chi) &=&\left\{\tilde{\varphi}_n^{[\ell]}(\chi)-r(\chi)f_n^{[\ell]}(\chi)-\Eb\left[\tilde{
\varphi}_n^{[\ell]}(\chi)-r(\chi)f_n^{[\ell]}(\chi)\right]\right\}\nonumber\\*
&&+\left\{ \Eb\left[\tilde{\varphi}_n^{[\ell]}(\chi)-r(\chi)f_n^{[\ell]}(\chi)\right]\right\}.
\end{eqnarray}
Noting that, from Lemma 3 in Amiri et al.  \cite{Amiri2}, we have $\Eb \left(f_n^{[\ell]}(\chi)\right)=M_1[1+o(1)]$ and it can be shown as the same lines of the proof of Lemma \ref{lem1}  that
\begin{eqnarray*}
f_n^{[\ell]}(\chi)-\Eb \left(f_n^{[\ell]}(\chi)\right)=O\left(\sqrt{\frac{\ln n}{n\phi(h_n)}}\right).
\end{eqnarray*}
Therefore, Theorem \ref{Cv_ps_rnl} follows from the combination of Lemmas \ref{Biais1} and \ref{lem1}, since from
Lemma \ref{Biais1}, if $ \displaystyle\lim\limits_{n\rightarrow+\infty}nh_n^{2}=0$, then
%\begin{eqnarray*}
$\varlimsup\limits_{n\rightarrow\infty}\left[\frac{n\phi(h_n)}{\ln~ n}\right]^{1/2}
\left\{ \Eb\left[\tilde{\varphi}_n^{[\ell]}(\chi)-r(\chi)f_n^{[\ell]}(\chi)\right]\right\}=0.$
%\end{eqnarray*}
%In addition, from Lemma 3 in  \cite{Amiri2}, we have $\Eb \left(f_n^{[\ell]}(\chi)\right)=M_1[1+o(1)]$ and it can be shown as the same lines of the proof of Lemma \ref{lem1}  that
%\begin{eqnarray*}
%f_n^{[\ell]}(\chi)-\Eb \left(f_n^{[\ell]}(\chi)\right)=O\left(\sqrt{\frac{\ln n}{n\phi(h_n)}}\right).
%\end{eqnarray*}  \hfill $\square$
%we have
%\begin{eqnarray*}
%U_n &\leq&  \frac{8(1+\kappa)}{\varepsilon\kappa}\frac{\|K\|_\infty }{f_1(\chi)B_{n,1-\ell}n\phi(h_n)}\gamma n^{1-p_0\rho}(\ln n)^{1/\mu}\\ & \leq & C n^{1-p_0\rho}(\ln n)^{1/\mu},
%\end{eqnarray*}
%so that $\sum_n U_n<\infty$ whenever $p_0>\frac{2}{\rho}$.
%Next observe that $\sum_n V_n<\infty$ whenever $\varepsilon>2(1+\kappa)V(\chi)$, and the conclusion follows from the application of Borel Cantelli's lemma.\\
%\subsubsection*{Proof of Lemma \ref{lem3}}

\subsection{Proof of Theorem \ref{mse}}
The mean square error of $r_n^{[\ell]}$ can be decomposed as follow:
%$$r_n^\ell-r=(r_n^\ell-r)\frac{\Eb f_n^\ell-f_n^\ell}{\Eb f_n^\ell}+\frac{r(\Eb f_n^\ell-f_n^\ell)}{\Eb f_n^\ell}+\frac{\varphi_n^\ell-r\Eb f_n^\ell}{\Eb f_n^\ell}$$
\begin{eqnarray*}
\Eb\left[\left(r_n^{[\ell]}({\chi})-r(\chi)\right)^2\right]&=&\frac{r^2(\chi)\text{Var}\left(f_n^{[\ell]}({\chi})\right)}{\Eb^2(f_n^{[\ell]}(\chi))}-\frac{2r(\chi)\cov\left(\phi_n^{[\ell]}(\chi),f_n^{[\ell]}(\chi)\right)}{\Eb^2(f_n^{[\ell]}(\chi))}\\
 & &+\frac{\text{Var}\left(\varphi_n^{[\ell]}({\chi})\right)}{\Eb^2(f_n^{[\ell]}(\chi))}+\frac{\Eb\varphi_n^\ell(\chi)-r(\chi)\Eb f_n^{[\ell]}(\chi)}{\Eb^2 f_n^{[\ell]}(\chi)}\\
 &&+o(h_n^2)+o\left(1/\left(n\phi(h_n)\right)\right).
 %& &+\frac{1}{\Eb^2(f_n^{[\ell]}(\chi))}\Eb\left(\left[r_n^{[\ell]2}(\chi)-r^2(\chi)\right]\left[f_n^{[\ell]}(\chi)-\Eb f_n^{[\ell]}(\chi)\right]^2\right)\\
% & & -\frac{2}{\Eb^2(f_n^{[\ell]}(\chi))}\Eb\left(\left[r_n^{[\ell]}(\chi)-r(\chi)\right]\left[\varphi_n^{[\ell]}(\chi)-\Eb\varphi_n^{[\ell]}(\chi)\right]\left[f_n^{[\ell]}(\chi)-\Eb f_n^{[\ell]}(\chi)\right]\right)\\ &:=&E_1-E_2+E_3+E_4-E_5.
\end{eqnarray*}
Theorem \ref{mse} follows from   Lemmas \ref{Biais1} - \ref{Var}.   \hfill $\square$
% we have to show the following claims:
%From Lemma \ref{Var}, our result will entirely proven if we show

\label{lastpage}
\end{document}